\documentclass[review]{elsarticle}
\usepackage{amssymb}
\usepackage{lscape}
\usepackage{amsthm}
\usepackage{graphicx}
\usepackage{caption2}
\usepackage{graphics}
\usepackage{lineno,hyperref,color,makeidx}
\modulolinenumbers[5]

\journal{Journal of \LaTeX\ Templates}

\begin{document}

\begin{frontmatter}

\title{Tails and probabilities for extreme outliers}
\tnotetext[mytitlenote]{Submitted to AIP Conference proceedings}
\author[mymainaddress,mycorrespondingauthor]{Pavlina Jordanova\corref{mycorrespondingauthor}}
\cortext[c1]{115 Universitetska str., 9712 Shumen, Bulgaria}
\ead{pavlina\_kj@abv.bg}

\author[mysecondaryaddress]{Monika P. Peteva}
\cortext[c2]{ 5 "James Bourchier" blvd., 1164 Sofia, Bulgaria.}

\address[mymainaddress]{Faculty of Mathematics and informatics, Shumen University, Bulgaria}
\address[mysecondaryaddress]{Faculty of Mathematics and Informatics, Sofia University, Bulgaria}

\begin{abstract}
The task of estimation of the tails of probability distributions having small samples seems to be still opened and almost unsolvable. The paper tries to make a step in filling this gap. In 2017 Jordanova et al. introduce six new characteristics of the heaviness of the tails of theoretical distributions. They rely on the probability to observe {\color{blue}mild or} extreme outliers. The main their advantage is that they always exist. This work presents some new properties of these characteristics. Using them six distribution sensitive estimators of the extremal index are defined. A brief simulation study compares their quality with the quality of Hill, t-Hill, Pickands and Deckers-Einmahl-de Haan estimators.
\end{abstract}

\begin{keyword}
\texttt{}\sep \LaTeX\sep Elsevier \sep template
\end{keyword}

\end{frontmatter}


\section{INTRODUCTION AND PRELIMINARIES}

One of the main tasks in Extreme value theory is estimation of extremal index. Given a huge sample it is solved by Hill \cite{hill1975simple}, t-Hill\cite{stehlik2010favorable, jordanova2012weak}, Pickands\cite{pickands1975statistical} and Deckers-Einmahl-de Haan\cite{dekkers1989moment} estimators. However their rates of convergence are fast only in case when the tail of the observed distribution is very close to Pareto one. The last makes difficult the task for estimating extremal index based on small samples. A good experience with this can be done when you try to estimate the tail index of the Hill-horror distribution. This distribution is discussed e.g. in Embrechts et al. (2013) \cite{embrechts2013modelling} or Resnick (2007) \cite{resnick2007heavy}. Therefore a preliminary classification of the tails of the distributions that can be used for preparing later on distribution sensitive estimators of the extremal index seems to be reasonable and very useful.
According to Klugman \cite{klugman2012loss} "The tail of a distribution ... is that part that reveals probabilities about large values". And now the question: "What does it mean "large values"?" arises. In order to clarify this concept we follow Tukey at al. \cite{tukey1977exploratory} and McGill et al. (1978) \cite{mcgill1978variations} approach. They define mild and extreme outliers and box-plots. The main statistics that they use are the quartiles of the empirical distribution and the interquartile range. In 2017 Jordanova et al. \cite{jordanova2017measuring} use their results and make classification of the probability distributions with respect to heaviness of their tails. They are based on the probability of the event to observe extreme outlier in a sample of independent observations. Analogously to the situations when we consider mean values and variances it is possible one distribution to belong to more than one distributional type with respect to this classification. However it shows us the most appropriate classes of distributions for fitting the corresponding distributional tail. In Section 2 new properties of these characteristics are obtained. The main their advantages are that they always exist and they are invariant with respect to increasing affine transformations. In Section 3 a new estimator of the extremal index is obtained and its properties are compared with the properties of  Hill \cite{hill1975simple}, t-Hill\cite{stehlik2010favorable, jordanova2012weak}, Pickands\cite{pickands1975statistical} and Deckers-Einmahl-de Haan\cite{dekkers1989moment} estimators. A beautiful summary of their properties could be found e.g. in Resnick et al. (2007) \cite{resnick2007heavy} or Embrechts et al. (2013) \cite{embrechts2013modelling}, and the references there in. The paper finishes with some conclusive remarks.

Along the paper $X_1, X_2, ..., X_n$ are independent identically distributed(i.i.d.) observations on a random variable (r.v.) $X$. Denote their cumulative distribution function (c.d.f.) by $F_X(x) = P(X \leq x)$, the theoretical p-quantiles by $F_X^\leftarrow(p) = inf\{x \in R: F_X(x) \geq p\}$, $p \in (0, 1]$, and the corresponding increasing order statistics by $X_{(1:n)} \leq X_{(2:n)} \leq ... \leq X_{(n:n)}$. There are many different definitions of the empirical p-quantiles $\hat{F}_X^\leftarrow (p)$. They can be found e.g. in Hyndman et al. (1996) \cite{hyndman1996sample}, Langford (2006)\cite{langford2006quartiles} or Parzen (1979) \cite{parzen1979nonparametric}. We use the following one
\begin{equation}\label{EmpiricaPQuantile}
\hat{F}_X^\leftarrow (p) = X_{([(n+1)p]:n)} + \{(n+1)p - [(n+1)p]\}\{X_{([(n+1)p]+1:n)} - X_{([(n+1)p]:n)}\}.
\end{equation}
Here $[a]$ means the integer part of $a$ and $\frac{1}{n+1} \leq p \leq \frac{n}{n+1}$. This definition entails $\hat{F}_X^\leftarrow (\frac{k}{n+1}) = X_{(k: n)}$ and the fact that the empirical quantile function is linearly interpolated between these points. This estimator is implemented in function $quantile$ in R (2018)\cite{R} as {\bf type = 6}. Arnold et al. (1992)\cite{arnold1992first}, Section 5.5 shows that $X_{([(n+1)p]:n)}$ is asymptotically unbiased estimator for $F^\leftarrow(p)$. According to \cite{embrechts2013modelling} if $k_n \to \infty$ and $k_n/n \to p \in (0, 1)$, for $n \to \infty$, then $X_{(k_n:n)} \to F_X^\leftarrow(p) $ almost sure. The last means that  sample quantiles are strongly consistent estimators of the theoretical quantiles $F^\leftarrow(p)$.  Arnold et al. (1992)\cite{arnold1992first} Th. 8.5.1. and Smirnov (1949) \cite{smirnov1949limit} find conditions for their asymptotic normality. Pancheva (1984) \cite{pancheva1155limit} is the first who describes the limiting probability laws for non-linearly normalized extreme order statistics. Pancheva and Gacovska (2014) \cite{pancheva2014asymptotic} investigate asymptotic behavior of central order statistics under monotone normalizations.  Their limit theorems propose further development of the results in this paper for different numbers of the central order statistics.
Recently Barakat et al. (2017) \cite{barakat2017new} model maxima under linear-power normalizations.

We will use these results for estimating the first  $Q_1(X) = F_X^\leftarrow (0.25)$ and the third $Q_3(X) = F_X^\leftarrow (0.75)$ quartile of $X$. They can be useful for estimating procedures based on relatively small samples because they are particular cases of the central order statistics and their rate of convergence seems to be faster than the rate of convergence of the extreme values.

\section{PROPERTIES OF $p_{eL}$, $p_{eR}$ AND $p_{e2}$ CHARACTERISTICS}

Here we consider the following three characteristics of extremely heavy left-, right- or two-sided tails of theoretical distributions introduced in Jordanova et al. (2017) \cite{jordanova2017measuring}:
\begin{enumerate}
  \item $p_{eL}(X) = P(X < Q_1(X) - 3IQR(X))$
  \item $p_{eR}(X) = P(X > Q_3(X) + 3IQR(X))$
  \item $p_{e2}(X) = p_{eL} + p_{eR},$
\end{enumerate}
where $IQR(X) = Q_3(X) - Q_1(X)$ is the inter quartile range of the theoretical distribution. It is clear that $\hat{Q}_1(X) - 3I\hat{Q}R(X)$ and  $\hat{Q}_3(X) + 3I\hat{Q}R(X)$ are weekly consistent L {\it estimators} correspondingly for $Q_1(X) - 3IQR(X)$ and $Q_3(X) + 3IQR(X)$. For general theory of L {\it estimators} see e.g. Arnold et al. (1992)\cite{arnold1992first}.
The next their properties show that they are invariant with respect to shifting to a constant or with respect to a product with a positive number. This makes them very prominent for differentiating heaviness of the tails of the distributions.

{\bf Theorem 1.} The characteristics $p_{eL}(X)$, $p_{eR}(X)$ and $p_{e2}(X)$ possess the following properties:
\begin{description}
  \item[a)] If $F_X(-x) = 1 - F_X(x)$, then $p_{eL}(X) = p_{eR}(-X)$.
  \item[b)] If $c = constant$, then  $p_{eL}(X) = p_{eL}(X + c)$, $p_{eR}(X) = p_{eR}(X + c)$, $p_{e2}(X) = p_{e2}(X + c)$.
  \item[c)] If the constant $c > 0$, then $p_{eL}(cX) = p_{eL}(X)$, $p_{eR}(cX) = p_{eR}(X)$, $p_{e2}(cX) = p_{e2}(X).$
  \item[d)] If $c < 0$, then $p_{eL}(cX) = p_{eR}(X), \quad p_{eR}(cX) = p_{eL}(X), \quad p_{e2}(cX) = p_{e2}(X).$
\end{description}

{\bf Sketch of the proof:} {\bf b)} is corollary of  the facts that
$$Q_1(X + c) = Q_1(X) + c, \quad Q_3(X + c) = Q_3(X) + c, \quad IQR(X + c) = IQR(X).$$

{\bf c)} follows by the equalities $Q_1(cX) = cQ_1(X)$, $Q_1(cX) = cQ_1(X)$ and $IQR(cX) = cIQR(X)$.

{\bf d)} Consider $c < 0$. Then $Q_1(cX) = cQ_3(X)$, $Q_3(cX) = cQ_1(X)$, $IQR(cX) = -cIQR(X)$.
\begin{eqnarray*}
  p_{eL}(cX) &=&  P(cX < Q_1(cX) - 3IQR(cX)) = P(cX < cQ_3(X) + 3cIQR(X)) \\
   &=& P(X > Q_3(X) + 3IQR(X)) = p_{eR}(X).
\end{eqnarray*}

\bigskip

In the next examples we will skip the cases when $p_{eL}(X) = 0$, $p_{eR}(X) = 0$ and $p_{e2}(X) = 0$ simultaneously. In this class of distributions fall e.g. Uniform distribution. The definitions of the distributions that we consider below could be found in many standard textbooks in probability theory.

\begin{figure}
\begin{minipage}[t]{0.5\linewidth}
\includegraphics[scale=.45]{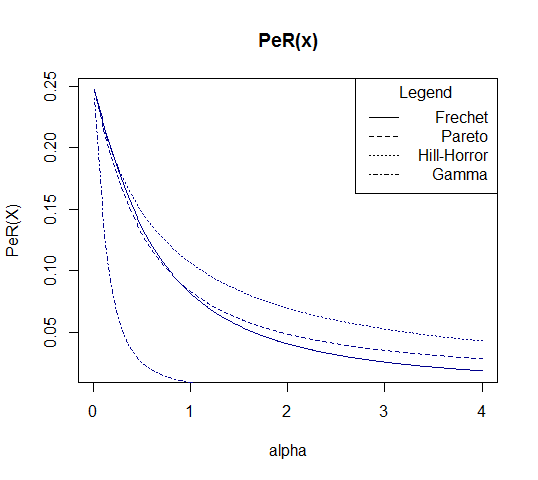}\vspace{-0.3cm}
    \end{minipage}
\begin{minipage}[t]{0.5\linewidth}
\includegraphics[scale=.45]{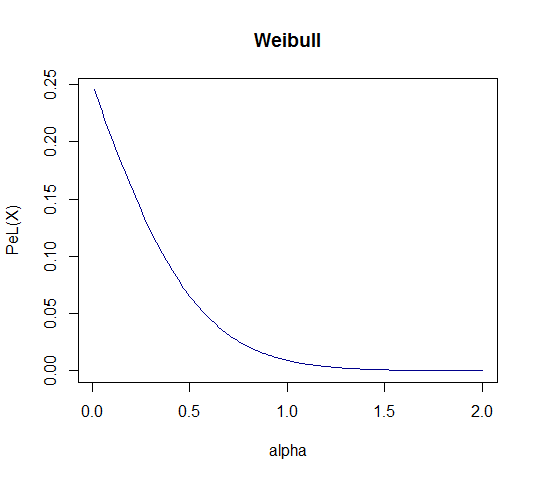}\vspace{-0.3cm}

\end{minipage}
\caption{The dependence of $p_{eR}(X)$ and $p_{eL}(X)$ from $\alpha$}
    \label{fig:PeLalpha}
\end{figure}

{\bf Example 1. Exponential distribution.} Let $\lambda > 0$, and $X \sim Exp(\lambda)$ with $EX = \frac{1}{\lambda}$. Then $F_X^\leftarrow(p) = -\frac{log(1 - p)}{\lambda}$, $Q_1(X) = \frac{log\, \frac{4}{3}}{\lambda}$, $Q_3(X) = \frac{log\, 4}{\lambda}$,  $IQR(X) = \frac{log\, 3}{\lambda}$,
$Q_1(X) - 3IQR(X) = \frac{log\, \frac{4}{3}}{\lambda} - 3\frac{log\, 3}{\lambda} = \frac{log\, 4 - 4 log\, 3}{\lambda} < 0$, $p_{eL}(X) = 0$,
$$Q_3(X) + 3IQR(X) = \frac{log\, 4}{\lambda} + 3\frac{log\, 3}{\lambda},$$
$$p_{e2}(X) = p_{eR}(X) = P(X > \frac{log\, 4}{\lambda} + 3\frac{log\, 3}{\lambda}) = \frac{1}{108} = 0.00925(925).$$

{\bf Example 2. Gamma distribution.} Assume $\alpha > 0$, $\beta > 0$, $X \sim Gamma(\alpha, \beta)$ with $EX = \frac{\alpha}{\beta}$.
Due to the fact that $\beta > 0$ is a scale parameter and the characteristics $p_{eL}(X)$ and $p_{eR}(X)$ are invariant with respect to a scale change,  without lost of generality $\beta = 1$. The probability for extreme left outliers $p_{eL}(X) = 0$. In order to obtain the quantile function of $X$ and to depict the dependence of $p_{eR}(X)$ on $\alpha$ we have used R software \cite{R}. The results are plotted on Figure \ref{fig:PeLalpha}, left. The conclusion that only in case $\alpha < 1$ we have $P_{eR} > 0$ corresponds to those made by Klugman et al. (2012) \cite{klugman2012loss} based on hazard rate function.

{\bf Example 3. Normal distribution.} Consider $\mu \in R$, $\sigma^2 > 0$ and $X \sim N(\mu, \sigma^2)$. Without lost of generality $\mu = 0$ and $\sigma^2 = 1$. Due to its symmetry $Q_1(X) = - Q_3(X) \approx -0.6745$ and $IQR(X) \approx 1.349$. Then
$Q_1(X) - 3IQR(X) = -4.7214 = -Q_3(X) - 3IQR(X)$. Therefore $p_{eL}(X) = p_{eR}(X) \approx 0.000001171.$

\begin{center}
\begin{table}
\begin{center}

\begin{tabular}{|c|c|c|c|c|c|}
  \hline
  n & 1 & 2 & 3 & 4 & 5   \\
  \hline
  $p_{eR}(X)= p_{eL}(X)$ & 0.0453 & 0.0146 & 0.0064 & 0.0033 & 0.0019 \\
    \hline
    \hline
     n &  6 & 7 & 8 & 9 & 10 \\
  \hline
  $p_{eR}(X)= p_{eL}(X)$ &  0.0012  & 0.0008 & 0.0006 & 0.0004 & 0.0003 \\
    \hline
\end{tabular}

\caption{The dependence of $p_{eR}(X)$ and $p_{eL}(X)$ from $n$. Here $X \sim t(n)$.}\label{tab:t}
\end{center}

\end{table}
\end{center}

{\bf Example 4. $t$-distribution.} Assume $n \in N$ and $X \sim t(n)$. Using R software \cite{R} one can obtain $Q_1(X) = - Q_3(X)$, $IQR(X)$ and $Q_1(X) - 3IQR(X) = -(Q_3(X) + 3IQR(X))$. The values of $p_{eR}(X)= p_{eL}(X)$ are presented in Table \ref{tab:t}.

{\bf Example 5. Pareto distribution.} Let $\alpha > 0$, $\delta > 0$ and
$$F_X(x) =\left\{\begin{array}{ccc}
                   0 & ,  & x < \delta \\
                   1 - \left(\frac{\delta}{x}\right)^\alpha & , & x \geq \delta
                 \end{array}
\right..$$

In this case $F_X^\leftarrow(p) = \frac{\delta}{\sqrt[\alpha]{1 - p}}$, $Q_1(X) = \frac{\delta\sqrt[\alpha]{4}}{\sqrt[\alpha]{3}}$, $Q_3(X) = \delta\sqrt[\alpha]{4}$,  $IQR(X) = \delta \sqrt[\alpha]{4}\left(1 - \frac{1}{\sqrt[\alpha]{3}}\right)$. Then
$$Q_1(X) - 3IQR(X) = \delta\sqrt[\alpha]{4}\left(\frac{4}{\sqrt[\alpha]{3}} - 3\right) \leq \delta,\,\, Q_3(X)  + 3IQR(X) = \delta\sqrt[\alpha]{4}\left(4 - \frac{3}{\sqrt[\alpha]{3}}\right),$$
$$ p_{eL}(X) = 0, \,\, p_{eR}(X) = \frac{3}{4(4\sqrt[\alpha]{3} - 3)^\alpha}.$$
The plot of the last characteristic with respect to $\alpha$ is presented on Figure \ref{fig:PeLalpha}, left.

{\bf Example 6. Fr$\acute{e}$chet distribution.}  Assume $\alpha > 0$, $\mu \in R$, $\sigma > 0$ and
$$F_X(x) =\left\{\begin{array}{ccc}
                   0 & ,  & x < \mu \\
                   exp\left\{-\left(\frac{x-\mu}{\sigma}\right)^{-\alpha}\right\} & , & x \geq \mu
                 \end{array}
\right..$$
It is well known that $F_X^\leftarrow(p) = \mu + \sigma (-log\, p)^{-1/\alpha}$, $Q_1(X) = \mu + \sigma (log\, 4)^{-1/\alpha}$, $Q_3(X) = \mu + \sigma (log\, 4 - log\, 3)^{-1/\alpha}$,  $IQR(X) = \sigma [(log\, 4 - log\, 3)^{-1/\alpha}-(log\, 4)^{-1/\alpha}]$. Then
$Q_1(X) - 3IQR(X) = \mu + \sigma [4(log\, 4)^{-1/\alpha} - 3 (log\, 4 - log\, 3)^{-1/\alpha}] ,$ \\
$Q_3(X)  + 3IQR(X) = \mu + \sigma\{4(log\, 4 - log\, 3)^{-1/\alpha} - 3(log\, 4)^{-1/\alpha}\},$
 $$p_{eL}(X) = \left\{\begin{array}{ccc}
                        0 & , & \alpha \in (0, 5.4662],  \\
                        exp\left\{-\left[4(log\,4)^{-1/\alpha} - 3(log\,\frac{4}{3})^{-1/\alpha}\right]^{-\alpha}\right\} & , & \alpha > 5.4662,
                      \end{array}
 \right.,$$
$$\frac{log\left(\frac{log\, 4}{log\left(\frac{4}{3}\right)}\right)}{log\left(\frac{4}{3}\right)} \approx 5.4662.$$
 $$p_{eR}(X) =  1 - exp\left\{-\left[4(log\,\frac{4}{3})^{-1/\alpha} - 3\,(log\,4)^{-1/\alpha}\right]^{-\alpha}\right\}.$$
The dependence of $p_{eR}(X)$ on $\alpha$ is depicted on Figure \ref{fig:PeLalpha}, left. It corresponds to the well known result that Fr$\acute{e}$chet's and Pareto's tails and very similar.

{\bf Example 7. Weibull negative distribution.}  Consider $\alpha > 0$, $\mu \in R$, $\sigma > 0$, and
$$F_X(x) =\left\{\begin{array}{ccc}
                   exp\left\{-\left(-\frac{x-\mu}{\sigma}\right)^{\alpha}\right\} & , & x < \mu\\
                   1 & ,  & x \geq \mu
                 \end{array}
\right..$$
 The corresponding quantile function is $F_X^\leftarrow(p) = \mu - \sigma (-log\, p)^{1/\alpha}$, $Q_1(X) = \mu - \sigma (log\, 4)^{1/\alpha}$, $Q_3(X) = \mu - \sigma (log\, 4 - log\, 3)^{1/\alpha}$,  and $IQR(X) = \sigma [(log\, 4)^{1/\alpha} - (log\, 4 - log\, 3)^{1/\alpha}]$. Then
$Q_1(X) - 3IQR(X) = \mu - \sigma [4(log\, 4)^{1/\alpha} - 3 (log\, 4 - log\, 3)^{1/\alpha}] ,$ $Q_3(X)  + 3IQR(X) = \mu + \sigma\{3(log\, 4)^{1/\alpha} - 4(log\, 4 - log\, 3)^{1/\alpha}\} > \mu,$
 $$p_{eL}(X) =   exp\left\{-\left[4(log\,4)^{1/\alpha} - 3(log\,\frac{4}{3})^{1/\alpha}\right]^{\alpha}\right\}, \quad p_{eR}(X) = 0. $$
Figure \ref{fig:PeLalpha}, right represents the dependence of $p_{eL}(X)$ on $\alpha$.

{\bf Example 8. Gumbell distribution.}  Let $\alpha > 0$, $\mu \in R$, $\sigma > 0$ and
$$F_X(x) = exp\left\{-exp\left[-\frac{x-\mu}{\gamma}\right]\right\},\quad  x \in R.$$
In this case $F_X^\leftarrow(p) = \mu - \gamma[log(-log\, p)]$, $Q_1(X) = \mu - \gamma[log(log\, 4)]$, $Q_3(X) = \mu - \gamma[log(log\, \frac{4}{3})]$,  and $IQR(X) = \gamma[log(log\, 4) - log(log\, \frac{4}{3})]$. Then
$$Q_1(X) - 3IQR(X) = \mu - \gamma\left[4\,log(log\,4) - 3log\left(log\frac{4}{3}\right)\right],$$
$$Q_3(X) + 3IQR(X) = \mu + \gamma\left[3\,log(log\,4) - 4log\left(log\frac{4}{3}\right)\right],$$
$$p_{eL}(X) = P\left(X < \mu-\gamma\left[4log(log\,4) - 3log\left(log\frac{4}{3}\right)\right]\right) = $$
$$ = exp\left\{-\frac{[log\,4]^4}{\left[log\left(\frac{4}{3}\right)\right]^3}\right\} \approx 4.264\times10^{-68}.$$
$$p_{eR}(X) = P\left(X > \mu + \gamma\left[3\,log(log\,4) - 4log\left(log\frac{4}{3}\right)\right]\right) =$$
 $$ = 1 - exp\left\{-\frac{[log\left(\frac{4}{3}\right)]^4}{\left[log\,4\right]^3}\right\} \approx 0.0026.$$
The last means that in the context of $p_{eR}(X)$ characteristics Exponential distribution has heavier right tail than the Gumbell one. Moreover it has  approximately three times higher chance for observing right extreme outliers.

{\bf Example 9. Hill-horror distribution.}  Assume $\alpha > 0$ and
$$F_X^\leftarrow(p) = \frac{-log(1 - p)}{\sqrt[\alpha]{1 - p}}, \quad p \in (0, 1).$$

Then $Q_1(X) = \frac{\sqrt[\alpha]{4}}{\sqrt[\alpha]{3}}log\,\frac{4}{3}$, $Q_3(X) = \sqrt[\alpha]{4}log\,4$,  $$IQR(X) = \sqrt[\alpha]{4}log\,4\left(1 - \frac{1}{\sqrt[\alpha]{3}}\right)+ \frac{\sqrt[\alpha]{4}}{\sqrt[\alpha]{3}}log\,3.$$
$$Q_1(X) - 3IQR(X) = \sqrt[\alpha]{4}\left(\frac{4}{\sqrt[\alpha]{3}}log\frac{4}{3} - 3log\, 4\right) < 0,$$
 $$Q_3(X)  + 3IQR(X) = \sqrt[\alpha]{4}\left(4log\, 4-\frac{3}{\sqrt[\alpha]{3}}log\frac{4}{3}\right) > 0.$$

Therefore $p_{eL}(X) = 0$. The dependence on the values of $p_{eR}(X)$ with respect to $\alpha$ is presented on Figure \ref{fig:PeLalpha}, left. We observe that within the considered types this distribution has "heaviest tail".






\section{THE EXTREMAL INDEX ESTIMATORS}

Suppose $X$ is a r.v. with c.d.f. $F_X$ with regularly varying tail. More precisely there exists $\alpha > 0$ such that for all $x > 0$,
\begin{equation}\label{RV}
  \lim_{t \to \infty} \frac{1 - F_X(xt)}{1 - F_X(t)} = x^{-\alpha}.
\end{equation}

 Jordanova et al. (2017) \cite{jordanova2017measuring} use $p_{eR}(X)$ characteristics and obtain five new distribution sensitive statistics for the parameter $\alpha$. Here we introduce one more estimator. It is based on the assumption that the observed r.v. has Hill-Horror distribution. Its rate of convergence is compared with the one of the corresponding  Hill \cite{hill1975simple}, t-Hill\cite{stehlik2010favorable, jordanova2012weak}, Pickands\cite{pickands1975statistical} and Deckers-Einmahl-de Haan\cite{dekkers1989moment} estimators.
Along the section $\hat{Q}_1(\vec{X}_n) = \hat{F}_X^\leftarrow (0.25)$ and $\hat{Q}_2(\vec{X}_n) = \hat{F}_X^\leftarrow (0.75)$ are correspondingly the first and the third empirical quartiles of the sample and $\hat{IQR}(\vec{X}_n) = \hat{Q}_3(\vec{X}_n) - Q_1(\vec{X}_n)$ is the empirical inter quartile range. The first group of two estimators that we consider are the most appropriate if the observed r.v. is Pareto distributed.

1.  Assume
$F_X(x) = 1  - x^{-\alpha}$,  $x > 0.$ Then $p_{eR}(X) = (Q_3(X) + 3(Q_3(X) - Q_1(X)))^{-\alpha}$, therefore
$$\alpha = -\frac{log\,\,p_{eR}(X)}{log\,\,(Q_3(X) + 3(Q_3(X) - Q_1(X)))}.$$

Having a sample of $n$ observations, if $\hat{Q}_1(\vec{X}_n) > 1$ the corresponding statistic is
$$\hat{\alpha}_{Par,n} = -\frac{log\,\,\hat{p}_{eR}(\vec{X}_n)}{log\,\,(\hat{Q}_3(\vec{X}_n) + 3(\hat{Q}_3(\vec{X}_n) - \hat{Q}_1(\vec{X}_n)))}.$$
where $\hat{p}_{eR}(\vec{X}_n)$ is the number of the extreme outlayers divided by the sample size $n$.

The quantiles are as useful as the cumulative distribution function. Quantile matching procedure seems to be well known. Its description could be seen e.g. in Klugman et al. (2012)\cite{klugman2012loss}. Analogously to the generalized method of moments we can make generalized quantile matching procedure. The second estimator is based on the fact that the fraction of the quartiles of Pareto distribution is
$\frac{Q_3(X)}{Q_1(X)} = \sqrt[\alpha]{3}.$
It is invariant with respect to a scale change, and given $\hat{Q}_1(\vec{X}_n)\not=\hat{Q}_3(\vec{X}_n)$ it has the form
$$\hat{\alpha}_{Par,Q} = \frac{log\,\,3}{log\,\,\hat{Q}_3(\vec{X}_n) - log\,\,\hat{Q}_1(\vec{X}_n)}.$$
Our empirical study shows that $\hat{\alpha}_{Par,Q}$ outperforms the other estimators discussed here in case of Pareto observed r.v.

2. The  estimators from the second group are the most appropriate for the case when the observed r.v. is Fr$\acute{e}$chet distributed.
The equality $F^\leftarrow(p_{eR}(X)) = Q_3(X) + 3(Q_3(X) - Q_1(X))$ leads us to the estimator
$$\hat{\alpha}_{Fr,n} = -\frac{log(-log(1-\hat{p}_{eR}(\vec{X}_n))}{log(\hat{Q}_3(\vec{X}_n) + 3(\hat{Q}_3(\vec{X}_n) - \hat{Q}_1(\vec{X}_n)))}.$$

In order to obtain the second estimator we consider the fraction
$$\frac{Q_3}{Q_1} = \frac{(-log\,\frac{3}{4})^{-1/\alpha}}{(-log\,\frac{1}{4})^{-1/\alpha}}.$$
Using the Generalized quantile matching procedure (see Klugman et al. (2012) \cite{klugman2012loss}) we express $\alpha$ and replace the theoretical quartiles with the corresponding empirical. Finally we obtain
$$\hat{\alpha}_{Fr,Q} = \frac{log\, (log\, 4) - log(log\,\frac{4}{3})}{log\,\,\hat{Q}_3(\vec{X}_n) - log\,\,\hat{Q}_1(\vec{X}_n)}.$$
Given small samples and Fr$\acute{e}$chet or Hill-Horror observed r.v. $\hat{\alpha}_{Fr,Q}$ seems to be very appropriate. It exceeds the quality of the other estimators discussed here. See Figure \ref{fig:4} and Figure \ref{fig:5}.

3. Suppose now that the observed r.v. has distribution which tail is close to those of the Hill-Horror distribution. In Embrechts et al. \cite{embrechts2013modelling} this distribution is defined via its quantile function
$F^\leftarrow(p) = (1 - p)^{-1/\alpha}(-log\,\,(1 - p))$, $p \in (0, 1).$

To the best knowledge of the authors the following estimator is new. It is obtained using the relation between $p_{eR}(X)$ characteristic of the Hill Horror distribution and $\alpha$. Given $p_{eR} \in (0, 1)$ and $Q_3(X) + 3(Q_3(X) - Q_1(X)) \not= -log\, p_{eR}$ we can express $\alpha$ and obtain
$$\hat{\alpha}_{HH,n} = \frac{log\, \hat{p}_{eR}(\vec{X}_n)}{log\,\frac{-log\,\hat{p}_{eR}(\vec{X}_n)}{\hat{Q}_3(\vec{X}_n) + 3(\hat{Q}_3(\vec{X}_n) - \hat{Q}_1(\vec{X}_n))}}.$$
In the next section we show that within the considered set of distributions $\hat{\alpha}_{HH,n}$ together with  $\hat{\alpha}_{Fr,Q}$ are the only appropriate estimators for $\alpha$ given small sample of observations on Hill-Horror distributed r.v. See Figure \ref{fig:5}.

It is easy to see that
$$\frac{Q_3(X)}{Q_1(X)} = \frac{3^{1/\alpha}(log\,\,4)}{log\,\,4-log\,\,3}.$$

Jordanova et al. (2017) replace the theoretical quartiles with the corresponding empirical. In this way using the generalized quantile matching procedure the authors define the following estimator
$$\hat{\alpha}_{HH,Q} = \frac{log\,\,3}{log\, \hat{Q}_3(\vec{X}_n) - log\, \hat{Q}_1(\vec{X}_n) + log\, log\, \frac{4}{3} - log \, log\, 4}.$$
The empirical results show that given a sample of observations on Hill-Horror distributed r.v. the rate of convergence of this estimator increases when $\alpha > 0$ decreases. However according to our observations  $\hat{\alpha}_{HH,Q}$ is too distribution sensitive and not robust.

\bigskip

\section{SIMULATION STUDY AND COMPARISON WITH ALTERNATIVE ESTIMATORS}

\begin{figure}
\begin{minipage}[t]{0.5\linewidth}
\includegraphics[scale=.49]{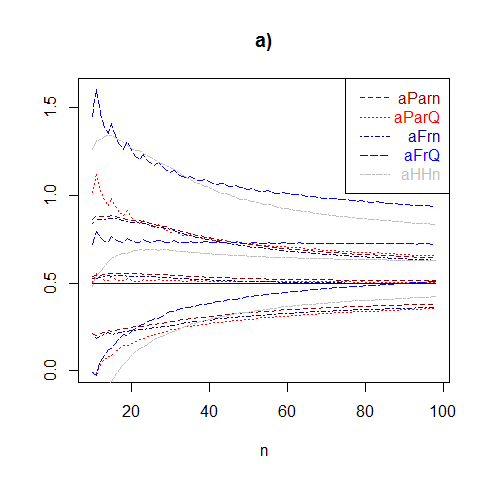}\vspace{-0.3cm}
    \end{minipage}
\begin{minipage}[t]{0.5\linewidth}
\includegraphics[scale=.49]{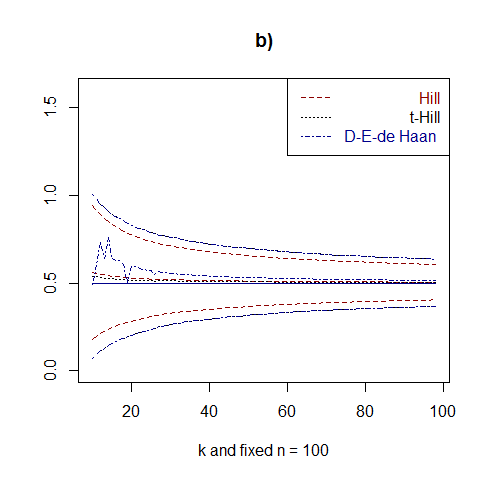}\vspace{-0.3cm}
\end{minipage}
\caption{Comparison between the rates of convergence of $\hat{\alpha}_{Par,n}$, $\hat{\alpha}_{Par,Q}$, $\hat{\alpha}_{Fr,n}$, $\hat{\alpha}_{Fr,Q}$, $\hat{\alpha}_{HH,n}$, (left) and  Hill, t-Hill and Deckers-Einmahl-de Haan, (right) estimators for Pareto(1, 0.5) observed r.v.}
    \label{fig:3}
\end{figure}

In this section we explore the behaviour of the considered estimators. Using the functions implemented in R (2018), \cite{R} we have simulated $m = 10000$ samples of $n = 100$ independent observations separately on r.v. with one of the following three probability laws: Pareto, Fr$\acute{e}$chet or Hill-Horror. Then for any fixed $n$ and for any fixed sample we have calculated $\hat{\alpha}_{Par,n}$, $\hat{\alpha}_{Par,Q}$, $\hat{\alpha}_{Fr,n}$, $\hat{\alpha}_{Fr,Q}$, $\hat{\alpha}_{HH,n}$  and $\hat{\alpha}_{HH,Q}$,  Hill \cite{hill1975simple}, t-Hill\cite{stehlik2010favorable, jordanova2012weak}, Pickands\cite{pickands1975statistical} and Deckers-Einmahl-de Haan\cite{dekkers1989moment} estimators. Finally we have averaged the corresponding values over the considered $n$. Because of Hill, t-Hill and Deckers-Einmahl-de Haan estimators depend not only of the sample size, but also from the number of order statistics that are included in their calculations in any of these three cases the results are plotted on separate figure. We have excluded the Pickands estimator from our plots because it turned out that the considered sample size $n = 100$ is not enough to observe its good properties.

On Figure \ref{fig:3}, a) is depicted the dependence of the values of $\hat{\alpha}_{Par,Q}$, $\hat{\alpha}_{Par,n}$, $\hat{\alpha}_{Fr,Q}$, $\hat{\alpha}_{Fr,n}$, $\hat{\alpha}_{HH,n}$ and their empirical 95\%  confidence intervals on the sample size. The plot of $\hat{\alpha}_{HH,Q}$ is  skipped because of it fluctuates too much. The plots of Hill, t-Hill, and Deckers-Einmahl-de Haan estimators together with their empirical 95\% confidence intervals are given on Figure \ref{fig:3}, b). Let us note that on the second figure the sample size is fixed. It is $n = 100$, and only the number of order statistics $k$ changes. Therefore this figure should be compared only with the points $n = 100$ on Figure \ref{fig:3}, a). We observe that $\hat{\alpha}_{Par,Q}$, $\hat{\alpha}_{Par,n}$ and $\hat{\alpha}_{Fr,n}$ have very similar behaviour to the well known estimators. Of course in this case, when the observed r.v. has exact Pareto distribution the Hill estimator outperforms the others.

\begin{figure}
\begin{minipage}[t]{0.5\linewidth}
\includegraphics[scale=.49]{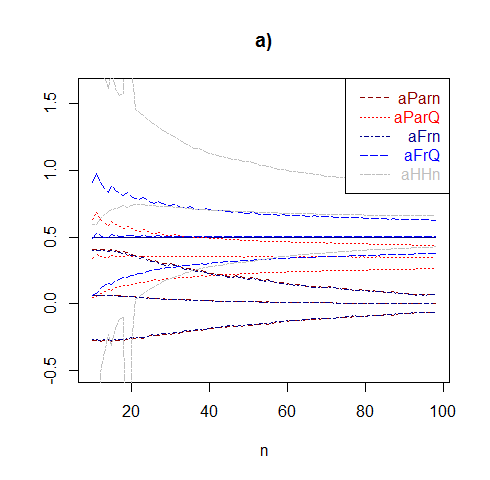}\vspace{-0.3cm}
    \end{minipage}
\begin{minipage}[t]{0.5\linewidth}
\includegraphics[scale=.49]{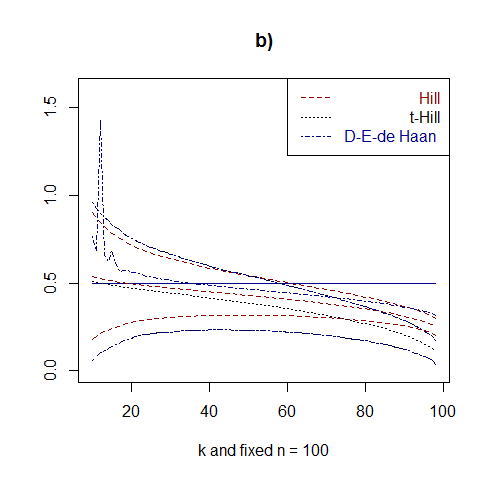}\vspace{-0.3cm}

\end{minipage}
\caption{Comparison between the rates of convergence of $\hat{\alpha}_{Par,n}$, $\hat{\alpha}_{Par,Q}$, $\hat{\alpha}_{Fr,n}$, $\hat{\alpha}_{Fr,Q}$, $\hat{\alpha}_{HH,n}$, $\hat{\alpha}_{HH,Q}$, (left) and  Hill, t-Hill and Deckers-Einmahl-de Haan, (right) estimators for Fr$\acute{e}$chet($\alpha = 0.5, \mu = 0, \sigma = 1$)  observed r.v.}
    \label{fig:4}
\end{figure}

If the observed r.v. is Fr$\acute{e}$chet($\alpha = 0.5, \mu = 0, \sigma = 1$) distributed, then the results from our simulation study, depicted on Figure \ref{fig:4}, a) show that only $\hat{\alpha}_{Fr,Q}$ estimator seems to be unbiased. The biggest advantage of this estimator is that in this case and for the considered sample sizes $n \leq 100$  Hill, t-Hill, Pickands and Deckers-Einmahl-de Haan estimators are not appropriate, because of their slower rate of convergence.

\begin{figure}
\begin{minipage}[t]{0.5\linewidth}
\includegraphics[scale=.44]{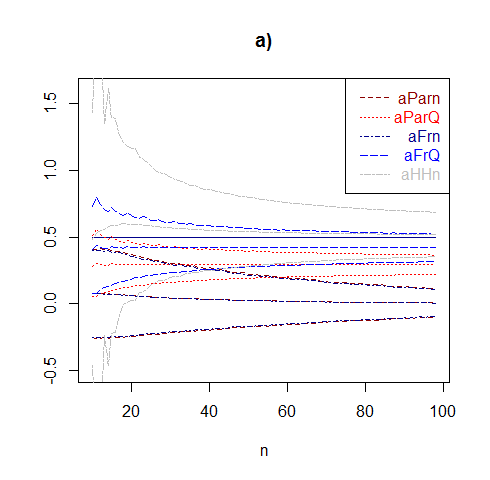}\vspace{-0.3cm}
    \end{minipage}
\begin{minipage}[t]{0.5\linewidth}
\includegraphics[scale=.44]{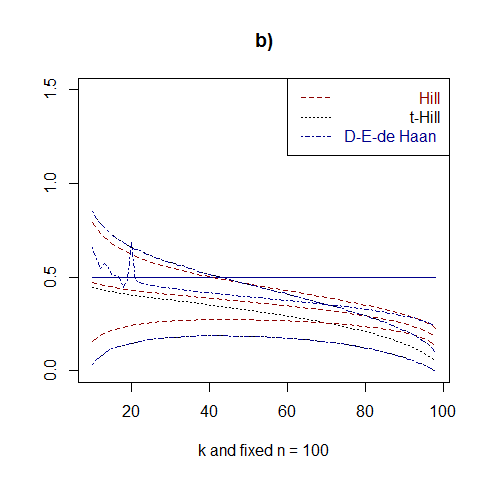}\vspace{-0.3cm}

\end{minipage}
\caption{Comparison between the rates of convergence of $\hat{\alpha}_{Par,n}$, $\hat{\alpha}_{Par,Q}$, $\hat{\alpha}_{Fr,n}$, $\hat{\alpha}_{Fr,Q}$, $\hat{\alpha}_{HH,n}$, $\hat{\alpha}_{HH,Q}$, (left) and  Hill, t-Hill and Deckers-Einmahl-de Haan, (right) estimators for Hill-Horror($\alpha = 0.5$)  observed r.v.}
    \label{fig:5}
\end{figure}

The case when the observed r.v. comes from Hill-Horror type is the most difficult for estimating. Here we have simulated such samples for $\alpha = 0.5$. Hill, t-Hill, Pickands and Deckers-Einmahl-de Haan statistics are not appropriate because the sample size $n = 100$ is too small. See e.g. Embrechts et al. (2013) \cite{embrechts2013modelling} and Figure  \ref{fig:5}, b). The plots on Figure \ref{fig:5}, a) show that only $\hat{\alpha}_{Fr,Q}$ and $\hat{\alpha}_{HH,n}$ estimators has relatively fast rate of convergence and seems to be appropriate in this case.

\bigskip

\section{CONCLUSIONS}

The introduced $p_{eL}$, $p_{eR}$ and $p_{e2}$ characteristics and their estimators are appropriate for usage in preliminary statistical analysis. They can help the practitioners to find the closest classes of probability laws to the distribution of the observed r.v. Within that family the tail index needs further estimation. That is when we fix the most appropriate parametric family the proposed estimators work well, but they are not appropriate in general non-parametric situations. For example if the observed r.v.  $X$ has Pareto distribution then it is well known that Hill estimator is the best one. Here we propose $\hat{\alpha}_{Par,n}$, $\hat{\alpha}_{Par,Q}$, $\hat{\alpha}_{Fr,n}$ estimators as its alternatives. In case when $X$ follows Fr$\acute{e}$chet type, then $\hat{\alpha}_{Fr,Q}$ has the best properties. If $X$ is close to Hill-Horror distribution $\hat{\alpha}_{HH,n}$ and $\hat{\alpha}_{Fr,Q}$ have fast rate of convergence and therefore they can be very useful for working with relatively small sample sizes. However the main disadvantage of all these estimators is that they are too distribution sensitive. The last means that their good properties disappear if the distributional type is not correctly determined. Here the characteristics of the heaviness of the tails of the distributions $p_{eL}(X)$, $p_{eR}(X)$, $p_{e2}(X)$, $p_{mL}(X)$, $p_{mR}(X)$ and $p_{m2}(X)$ can be useful.

\section*{ACKNOWLEDGMENTS}
The authors are grateful to the bilateral projects Bulgaria - Austria, 2016-2019, Feasible statistical modelling for extremes in ecology and finance, Contract number 01/8, 23/08/2017, the project RD-08-125/06.02.2018 from the Scientific Research Fund in Shumen University, and the project 80-10-222/04.05.2018 from the Scientific Funds of Sofia University.

\nocite{*}
\bibliographystyle{aipnum-cp}%

\end{document}